\documentclass[12pt]{amsart}
\usepackage{amsmath}
\usepackage{amsxtra}
\usepackage{amscd}
\usepackage{amsthm}
\usepackage{amsfonts}
\usepackage{amssymb}
\usepackage{eucal}
\usepackage{mathrsfs}
\usepackage{hyperref}

\def\C{{\mathbb C}}

\def\F{{\mathbb F}}

\def\P{{\mathbb P}}

\def\Z{{\mathbb Z}}

\newtheorem{theorem}{Theorem}[section]
\newtheorem{definition}{Definition}[section]
\newtheorem{lemma}[theorem]{Lemma}

\newtheorem{corollary}[theorem]{Corollary}

\title{Noncomplete embeddings of Rational Surfaces}
\author{Euisung Park }
\address {Euisung Park : School of Mathematics, Korea Institute for Advanced Study,
270-43 Cheongryangri-dong, Dongdaemun-gu, Seoul 130-012, Republic
of Korea,} \email{puserdos@kias.re.kr}
\thanks{The author was supported by Korea Research Foundation Grant (KRF-2002-070-C00003).}

\begin{document}

\begin{abstract}
In this paper, we study the Castelnuovo-Mumford regularity of
nonlinearly normal embedding of rational surfaces. Let $X$ be a
rational surface and let $L \in \mbox{Pic}X$ be a very ample line
bundle. For a very ample subsystem $V \subset H^0 (X,L)$ of
codimension $t \geq 1$, if $X \hookrightarrow \P (V)$ satisfies
Property $N^S _1$, then $\mbox{Reg} (X) \leq t+2$\cite{KP}. Thus
we investigate Property $N^S _1$ of noncomplete linear systems on
X. And our main result is about a condition of the position of $V$
in $H^0 (X,L)$ such that $X \hookrightarrow \P (V)$ satisfies
Property $N^S _1$. Indeed it is related to the geometry of a
smooth rational curve of $X$. Also we apply our result to $\P^2$
and Hirzebruch surfaces.
\end{abstract}

\maketitle \tableofcontents

\section{Introduction}
\noindent In this paper we study the Castelnuovo-Mumford
regularity of nonlinearly normal embedding of rational surfaces.
Let $X$ be a smooth projective variety. Each embedding of $X$
corresponds to a very ample line bundle $L \in \mbox{Pic}X$ and a
very ample subsystem $V \subset H^0 (X,L)$. Put $\mbox{dim}_{\C} V
=r+1$ and let $X \hookrightarrow \P^r$ be the embedding defined by
$V$. To understand this embedding, the classical problem is:\\

\begin{enumerate}
\item[$(*)$] Let $d \geq 1$ be a fixed integer.
\begin{enumerate}
\item[(1)] ($d$-normality) Determine whether the map $S^d H^0
(\P^r,\mathcal{O}_{\P^r}(d)) \rightarrow H^0 (X,L^d)$ is
surjective or not.
\item[(2)] Determine whether the degrees of all
minimal generators of the homogeneous ideal of $X$ is bounded by
$d$.\\
\end{enumerate}
\end{enumerate}

\noindent Nowadays the problem $(*)$ has been reformulated as
finding the upper bound of the \textit{Castelnuovo-Mumford}
regularity of the ideal sheaf $\mathcal{I}_X$ of $X
\hookrightarrow \P^r (V)$. We say that $X \hookrightarrow \P (V)$
is $m$-regular if $H^i (\P^r , \mathcal{I}_X (m-i))=0$ for every
$i \geq 1$. It is well-known that if $X \hookrightarrow \P^r$ is
$m$-regular, then it is $k$-normal for all $k \geq m-1$, $I_X$ is
generated by forms of degree $\leq m$ and hence there is no
$(m+1)$-secant line to $X$(Lecture 14 in \cite{M}). Recently
Sijong Kwak and the author\cite{KP} discovered a relation between
the Castelnuovo-Mumford regularity of $X \hookrightarrow \P^r$ and
the minimal free resolution of $E = \oplus_{\ell \in \Z} H^0
(X,L^{\ell})$ as a graded $S = Sym^{\bullet} V$. For precise
statements, we introduce the following:

\begin{definition}
Assume that $E$ has a minimal free resolution of the form
\begin{equation*}
\cdots \rightarrow \oplus_j S ^{k_{i,j}}(-i-j) \rightarrow \cdots
\rightarrow \oplus_j S^{k_{0,j}}(-j) \rightarrow E \rightarrow 0.
\end{equation*}
For a nonnegative integer $p$, $~X \subset \P^r$ is said to
satisfy Property $N^{S}_p$ if $k_{i,j} = 0$ for $0 \leq i \leq p$
and  $j \geq 2$.
\end{definition}

\noindent Thus Property $N^S _p$ means that $E$ admits a minimal
free resolution of the form
\begin{eqnarray*}
\cdots \rightarrow S^{k_{p,1}}(-p-1) \rightarrow \cdots
\rightarrow S^{k_{2,1}}(-3) \rightarrow  S^{k_{1,1}}(-2)
\rightarrow S \oplus S^{k_{0,1}} (-1) \rightarrow E \rightarrow 0.
\end{eqnarray*}

\noindent In this case it is obvious that $k_{0,1} = h^0 (X,L)
-r-1$. The importance of Property $N^S _p$ comes from the fact
that if $X \hookrightarrow \P^r$ satisfies Property $N^S _1$, then
$k$-normality holds for all $k \geq k_{0,1}+1$, the homogeneous
ideal is generated by forms of degree $\leq k_{0,1}+2$ and
$\mbox{max} \{k_{0,1}+2, m \}$-regularity holds where $m$ is the
regularity of $\mathcal{O}_X$ with respect to $L$.(Theorem 1.1 in
\cite{KP}). Also the following is proved:

\begin{theorem}[Theorem 1.2 in
\cite{KP}]\label{thm:effect} Let $X$ be a smooth projective
variety and $L \in \mbox{Pic}X$ a very ample line bundle such that
$H^1 (X,L^j)=0$ for all $j\geq2$. If $L$ satisfies Property $N_p$,
then for every subsystem $V \subset H^0 (X,L)$ of codimension
$\leq p-1$, $X \hookrightarrow \P (V)$ satisfies Property
$N^{S}_1$.
\end{theorem}

\noindent However if $L$ fails to satisfy Property $N_{p+1}$ and
$V \subset H^0 (X,L)$ is a very ample subsystem of codimension
$p$, then $X \hookrightarrow \P (V)$ may not satisfy Property $N^S
_1$. See Theorem 5.1 in \cite{KP}. Therefore it seems natural to
find a condition of $V$ such that $X \hookrightarrow \P (V)$
satisfies Property $N^S _1$. And this paper is devoted to consider
this problem when $X$ is a smooth projective rational surface. The
following is our main result:

\begin{theorem}\label{thm:observationrational}
Let $X$ be a smooth projective rational surface and $L \in
\mbox{Pic}X$ a very ample line bundle such that $H^1 (X,L^j)=0$
for all $j \geq 1$. Let $V \subset H^0 (X,L)$ be a very ample
subsystem of codimension $t \geq 1$. Assume that there exists a
smooth rational curve $C \subset X$ satisfying the followings:
\begin{enumerate}
\item[$(i)$] The self-intersection number $C^2 =m$ is $\geq 2$.
\item[$(ii)$] The image of the map $V \rightarrow H^0 (C,L \otimes
\mathcal{O}_C)$ has codimension $\leq m-2$.
\end{enumerate}
Then $X \hookrightarrow \P (V)$ satisfies Property $N^S _1$ and
thus $(t+2)$-regularity holds.
\end{theorem}

\begin{corollary}\label{cor:generalsubsystem}
Let $X$ be a smooth projective rational surface and $L \in
\mbox{Pic}X$ a very ample line bundle such that $H^1 (X,L^j)=0$
for all $j \geq 1$. Assume that there exists a smooth rational
curve $C \subset X$ with the self-intersection number $C^2 =m \geq
2$. Then Property $N^S _1$ holds for general very ample subsystem
$V \subset H^0 (X,L)$ of codimension $\leq h^0 (X,L\otimes
\mathcal{O}_X (-C)) +m-2$.
\end{corollary}

\noindent {\bf Remark 1.1.} Let $X$ be a smooth surface with $p_g
= q =0$ and $L \in \mbox{Pic}X$ a very ample line bundle. If $-K_X
\cdot L > 0$, then $H^1 (X,L^j)=0$ for all $j \geq 1$. See Lemma
1.2 in \cite{GP}. In particular, this is the case when
$X$ is an anticanonical rational surface.\qed\\

\noindent {\bf Remark 1.2.} For the Castelnuovo-Mumford regularity
of surfaces, there is a very general result. Indeed , R.
Lazarsfeld\cite{L} proved that for a smooth projective surface $X
\hookrightarrow \P^r$ of degree
$d$ , $\mbox{Reg}(X) \leq d-r+3$. \qed \\

\noindent To apply Theorem \ref{thm:observationrational} to a
given rational surface $X$ and a very ample line bundle on it, we
need to choose an appropriate smooth rational curve on $X$. When
$X =\P^2$, let $C$ be a smooth quadric in $|\mathcal{O}_{\P^2}
(2)|$. Then for any $a \geq 2$, Theorem
\ref{thm:observationrational} can be applied to
$L=\mathcal{O}_{\P^2} (a)$. When $X = \F _e$, for any $e \geq 0$,
is the Hirzebruch surface defined by $\mathcal{O}_{\P^1} \oplus
\mathcal{O}_{\P^1} (-e)$, let $C_0$ be the minimal section and a
fiber $f$ by following the notation and terminology of R.
Hartshorne's book \cite{H}, V $\S 2$. Then $C_0 + (e+1)f$ is a
very ample line bundle and a general section $C \in |C_0 +
(e+1)f|$ is a smooth rational curve with $C^2 = e+2$. Note that
for every very ample line bundle $L \in \mbox{Pic}X$, $H^1
(X,L^j)=0$ for all $j \geq 1$. Therefore Theorem
\ref{thm:observationrational} can be applied to every very ample
line bundle on $X$. Finally for a rational surface $X$ with $K_X
^2 \leq 7$, there is a birational morphism $X \rightarrow \F_e$
for some $e \geq 0$. Therefore we can choose $C$ as the preimage
of a general section in $|C_0 + (e+1)f|$. In particular, we can
apply Theorem \ref{thm:observationrational} to every very ample
line bundle on anticanonical rational surfaces by Remark 1.1.

The organization of this paper is as follows. In $\S 2$, we review
some cohomological nature of Property $N^S _p$. $\S 3$ is devoted
to prove Theorem \ref{thm:observationrational} and Corollary
\ref{cor:generalsubsystem}. In $\S 4$, we apply our results to
$\P^2$ and Hirzebruch surfaces.\\

\section{The Property $N^S _p$}
\noindent Let $X$ be a projective variety and let $L \in
\mbox{Pic}X$ be a very ample line bundle. For the complete
embedding $X \hookrightarrow \P H^0 (X,L)$, Property $N_p$
(introduced by M. Green and R. Lazarsfeld\cite{GL}) has been
studied by several authors. In \cite{KP}, this notion is
generalized to Property $N^S _p$ for arbitrary very ample linear
systems. Let $V \subset H^0 (X,L)$ be a very ample subsystem of
codimension $t$ which defines the embedding $X \hookrightarrow \P
(V)$. Let $S$ be the homogeneous coordinate ring of $\P(V)$. Then
$L$ defines the graded $S$-algebra $E = \oplus_{\ell \in \Z} H^0
(X, L^{\ell})$. Many algebraic geometers studied the syzygies of
$E$ as a graded $\mbox{Sym}^\bullet H^0(X,L)$-module, i.e. when
$L$ satisfies Property $N_p$. This can be generalized to
noncomplete embeddings as follows:

Note that for complete embedding Property $N^{S}_p$ is equal to
Property $N_p$. Although the homogenous coordinate ring of $X
\subset \P(V)$ is a proper submodule of $E$ for $t=k_{0,1} \geq
1$, there is the following influence of Property $N^{S}_1$ on
higher normality and defining equations:

\begin{theorem}[Theorem 1.1, \cite{KP}]\label{thm:regularity}
Let $X$ be a smooth complex projective variety and let $L \in
\mbox{Pic}X$ be a very ample line bundle. For an embedding $ X
\hookrightarrow \P(V)$ given by a very ample subsystem $V \subset
H^0 (X,L)~~ $ of codimension $t$, assume that Property $N^{S}_1$
holds. Then $X \hookrightarrow \P(V)$ is
\begin{enumerate}
\item[$(1)$] $k$-normal for all $k \geq t+1$, \item[$(2)$]
max$\{m+1,t+2\}$-regular and, \item[$(3)$] the homogeneous ideal
of $X \hookrightarrow \P(V)$ is generated by forms of degree $\leq
t+2$.
\end{enumerate}
where $m$ is the regularity of $\mathcal{O}_X$ with respect to
$L$.
\end{theorem}

\noindent  Property $N^S _p$ is described cohomologically as
follows:

\begin{lemma}\label{lem:criterion}
Assume that $H^1 (X,L^j)=0$ for all $j \geq 1$. For a very ample
subsystem $V \subset H^0 (X,L)$, define the vector bundle
$\mathcal{M}_V$ by the following exact sequence
\begin{equation*}
0 \rightarrow \mathcal{M}_V \rightarrow V \otimes \mathcal{O}_X
\rightarrow L \rightarrow 0.
\end{equation*}
For $p \geq 0$, the embedding $X \hookrightarrow \P(V)$ satisfies
Property $N^S _p$ if and only if
\begin{equation*}
H^1 (X, \wedge^{i+1} \mathcal{M}_V \otimes L^j) =0~~~~\mbox{for $0
\leq i \leq p$ and all $j \geq 1$.}
\end{equation*}
\end{lemma}

\begin{proof}
See Lemma 3.2 in \cite{KP}.
\end{proof}

\noindent Also for surfaces that will concern us, we have the
following more simple criterion:

\begin{lemma}\label{lem:simplecriterion}
Let $X$ be a smooth projective surface with $p_g = q =0$ and
assume that $H^1 (X,L^j)=0$ for all $j \geq 1$. Then $X
\hookrightarrow \P(V)$ satisfies Property $N^S _p$ if and only if
\begin{equation*}
H^2 (X, \wedge^{i+1} \mathcal{M}_V )=0~~~~\mbox{for $1 \leq i \leq
p+1$.}
\end{equation*}
\end{lemma}

\begin{proof}
From the exact sequence
\begin{equation*}
0 \rightarrow \wedge^{i+1} \mathcal{M}_V  \rightarrow \wedge^{i+1}
V \otimes \mathcal{O}_X \rightarrow \wedge^i \mathcal{M}_V \otimes
L \rightarrow 0,
\end{equation*}
we know that $H^1 (X,\wedge^i \mathcal{M}_V \otimes L^{\ell +1})
\cong H^2 (X,\wedge^{i+1} \mathcal{M}_V \otimes L^{\ell})$ for all
$\ell \geq 0$. Also from the exact sequence
\begin{equation*}
0 \rightarrow \wedge^{i+2} \mathcal{M}_V  \rightarrow \wedge^{i+2}
V \otimes \mathcal{O}_X \rightarrow \wedge^{i+1} \mathcal{M}_V
\otimes L \rightarrow 0,
\end{equation*}
$H^2 (X,\wedge^{i+1} \mathcal{M}_V \otimes L^{\ell})=0$ for all
$\ell \geq 1$.
\end{proof}

\section{Proof of The Main Theorem}
\noindent This section is devoted to show Theorem
\ref{thm:observationrational}. Indeed the statement holds for
smooth projective surfaces with $p_g = q =0$. Thus we consider the
following situation throughout
this section:\\

\begin{enumerate}
\item[$(*)$] Let $X$ be a smooth projective surface with $p_g = q
=0$ and assume that there is a smooth rational curve $C \subset X$
such that the self-intersection number $C^2 :=m \geq 2$. Let $L
\in \mbox{Pic}X$ be a very ample line bundle such that $H^1
(X,L^j)=0$ for all $j \geq 1$. For the complete embedding $X
\hookrightarrow \P H^0 (X,L)$, assume that $C$ is embedded in the
linear subspace $\Lambda = \P (W) \subset \P H^0 (X,L)$ where $W
\subset H^0 (C,L \otimes \mathcal{O}_C)$ is the image of the
natural map $H^0 (X,L) \rightarrow H^0 (C,L \otimes
\mathcal{O}_C)$.\\
\end{enumerate}

\noindent Under the situation $(*)$, the aim of this section is to
prove the following:

\begin{theorem}\label{thm:observation}
For a very ample subsystem $V \subset H^0 (X,L)$ of codimension
$t$, consider the natural map $V \rightarrow H^0 (C,L \otimes
\mathcal{O}_C)$. If the image of $V$ has codimension $\leq m-2$,
then $X \hookrightarrow \P (V)$ satisfies Property $N^S _1$.
\end{theorem}

\noindent In order to prove these two theorems, we first show a
criterion for Property $N^S _1$.

\begin{lemma}
Let $X$ be a smooth projective surface with $p_g = q =0$. Let $C
\subset X$ be a smooth curve and put $A:= \mathcal{O}_X (C) \in
\mbox{Pic}X$. For a very ample line bundle $L \in \mbox{Pic}X$,
assume that
\begin{equation*}
H^1 (X,L^j)=0~~~~\mbox{for all $j \geq 1$}
\end{equation*}
For a very ample subsystem $V \subset H^0 (X,L)$, consider the
short exact sequence
\begin{equation*}
0 \rightarrow \mathcal{M}_V \rightarrow V \otimes \mathcal{O}_X
\rightarrow L \rightarrow 0.
\end{equation*}
If $H^1 (C,\wedge^{i+2} \mathcal{M}_V \otimes A^{\ell +1} \otimes
\mathcal{O}_C)=0$ for $0 \leq i \leq p$ and all $\ell \geq 0$,
then $X \hookrightarrow \P (V)$ satisfies Property $N^S _p$.
\end{lemma}

\begin{proof}
By Lemma \ref{lem:simplecriterion}, it suffices to check that
\begin{equation*}
H^2 (X, \wedge^{i+1} \mathcal{M}_V )=0~~~~\mbox{for $1 \leq i \leq
p+1$.}
\end{equation*}
From the exact sequence
\begin{equation*}
0 \rightarrow \wedge^{i+2} \mathcal{M}_V  \rightarrow \wedge^{i+2}
\mathcal{M}_V \otimes A \rightarrow \wedge^{i+2} \mathcal{M}_V
\otimes A \otimes \mathcal{O}_C \rightarrow 0,
\end{equation*}
we obtain
\begin{equation*}
H^1 (C,\wedge^{i+2} \mathcal{M}_V \otimes A^{\ell +1} \otimes
\mathcal{O}_C) \rightarrow H^2 (X,\wedge^{i+2} \mathcal{M}_V
\otimes A ^{\ell}) \rightarrow H^2(X,\wedge^{i+2} \mathcal{M}_V
\otimes A^{\ell +1}) \rightarrow 0.
\end{equation*}
Therefore $H^2 (X,\wedge^{i+2} \mathcal{M}_V \otimes A ^{\ell})
\cong H^2(X,\wedge^{i+2} \mathcal{M}_V \otimes A^{\ell +1})$ for
all $\ell \geq 0$. Also as proved in the following Lemma
\ref{lem:vanishing}, $H^2 (X,\wedge^{i+2} \mathcal{M}_V \otimes A
^{\ell})=0$ for all sufficiently large $\ell$. We thus have the
desired vanishing.
\end{proof}

\begin{lemma}\label{lem:vanishing}
Let $X$ be a smooth projective surface. Let $L \in \mbox{Pic}X$ be
an ample line bundle and let $A \in \mbox{Pic}X$ be a nontrivial
effective line bundle.\\
$(1)$ $H^2 (X,tA-L) = 0$ for all sufficiently large $t$.\\
$(2)$ Assume that $L$ is very ample and let $V \subset H^0 (X,L)$
be a base point free subsystem.

Consider the short exact sequence
            \begin{equation*}
             0 \rightarrow \mathcal{M}_V \rightarrow V \otimes \mathcal{O}_X \rightarrow
             L \rightarrow 0.
             \end{equation*}

Then for each $q \geq 1$, $H^2 (X,\wedge^q \mathcal{M}_V \otimes A
^t)=0$ for all sufficiently large $t$.
\end{lemma}

\begin{proof}
$(1)$ By Serre duality, $H^2 (X,tA-L) \cong H^0 (X,K_X +L
-tA)^\vee$. Note that since $L$ is ample and $A$ is nontrivial
effective, $L \cdot A >0$. Thus for all sufficiently large $t$,
\begin{equation*}
L \cdot (K_X + L - tA)=-t L \cdot A + L \cdot (L+K_X) < 0
\end{equation*}
and hence $H^0 (X,K_X +L-tA)=0$.\\
$(2)$ Consider the short exact sequence
\begin{equation*}
0 \rightarrow \wedge^{q+1} \mathcal{M}_V \otimes L^{-1} \otimes
B^t \rightarrow \wedge^{q+1}V \otimes L^{-1} \otimes A^t
\rightarrow \wedge^q \mathcal{M}_V \otimes A^t \rightarrow 0.
\end{equation*}
The assertion follows immediately from $(1)$.
\end{proof}

\noindent {\bf Proof of Theorem \ref{thm:observation}.} By Lemma
\ref{lem:simplecriterion}, it  suffices to show that $H^1 (C,
\wedge^{i+2} \mathcal{M}_V \otimes A^{\ell +1} \otimes
\mathcal{O}_C )=0$ for $i=0,1$ and all $\ell \geq 0$. Note that
$A|_C = \mathcal{O}_{\P^1} (m)$. Let $V' \subset H^0 (C,L \otimes
\mathcal{O}_C)$ be the image of $V \rightarrow H^0 (C,L \otimes
\mathcal{O}_C)$ and consider the short exact sequence
\begin{equation*}
0 \rightarrow \mathcal{M}_{V'} \rightarrow V' \otimes
\mathcal{O}_C \rightarrow L\otimes \mathcal{O}_C \rightarrow 0.
\end{equation*}
Since $\mbox{codim} (V',H^0 (C,L \otimes \mathcal{O}_C)) \leq
m-2$,
\begin{equation*}
H^1 (\P^1,\wedge^i \mathcal{M}_V \otimes \mathcal{O}_{C} (j))=H^1
(C,\wedge^i \mathcal{M}_V \otimes \mathcal{O}_{\P^1} (jm))=0
\end{equation*}
for $i=1,2,3$ if $1 \geq m$ by the following Lemma
\ref{lem:rationalcurve}. The second statement comes immediately
from Theorem \ref{thm:regularity} since $\mathcal{O}_X$ is
$2$-regular with respect to $L$.\qed\\

\begin{lemma}\label{lem:rationalcurve}
Let $V \subset H^0 (\P^1,\mathcal{O}_{\P^1} (d))$ be a base point
free subsystem of codimension $t$. Consider the short exact
sequence
\begin{equation*}
0 \rightarrow \mathcal{M}_V \rightarrow V \otimes
\mathcal{O}_{\P^1} \rightarrow \mathcal{O}_{\P^1} (d) \rightarrow
0.
\end{equation*}
Then $H^1 (\P^1,\wedge^i \mathcal{M}_V \otimes \mathcal{O}_{\P^1}
(j))=0$ for $i=1,2,3$ if $j \geq t+2$.
\end{lemma}

\begin{proof}
Consider the following commutative diagram:\\
\begin{equation*}
\begin{CD}
              & 0                  &             & 0                                &             &           & \\
              & \downarrow         &             & \downarrow                       &             &           & \\
0 \rightarrow & \mathcal{M}_{V} & \rightarrow & V \otimes \mathcal{O}_{\P^1}        & \rightarrow & \mathcal{O}_{\P^1} (d)& \rightarrow 0 \\
              & \downarrow         &             & \downarrow                       &             & \parallel & \\
0 \rightarrow & \mathcal{M}_{d}   & \rightarrow & H^0 (C,\mathcal{O}_{\P^1} (d)) \otimes \mathcal{O}_{\P^1} & \rightarrow & \mathcal{O}_{\P^1} (d)& \rightarrow 0 \\
              & \downarrow         &             & \downarrow                       &             &           & \\
              & \mathcal{O}_{\P^1} ^t   &   =    & \mathcal{O}_{\P^1} ^t            &             &           & \\
              & \downarrow         &             & \downarrow                       &             &           & \\
              & 0                  &             & 0                                &             &           &.
\end{CD}
\end{equation*}\\
Since every vector bundle over $\P^1$ is totally decomposable,
$\mathcal{M}_d = \oplus^d \mathcal{O}_{\P^1} (-1)$ and
\begin{equation*}
\mathcal{M}_V = \oplus_{i=1} ^{d-t} \mathcal{O}_{\P^1}
(-a_i)~~\mbox{where $1\leq a_1 \leq \cdots \leq a_{d-t}$ and $a_1
+ \cdots + a_{d-t} = d$.}
\end{equation*}
Since $a_{i_1} + a_{i_2} +a_{i_3} \leq t+3$ for distinct
$i_1,i_2,i_3$, $H^1 (\P^1,\wedge^3 \mathcal{M}_V \otimes
\mathcal{O}_{\P^1} (j))=0$ if $j \geq t+2$. By the same way, one
can check the vanishing of $H^1 (\P^1,\wedge^2 \mathcal{M}_V
\otimes \mathcal{O}_{\P^1} (j))$ and $H^1 (\P^1,\mathcal{M}_V
\otimes \mathcal{O}_{\P^1} (j))$ for $j \geq t+2$.
\end{proof}

\noindent {\bf Proof of Theorem \ref{thm:observationrational}.}
Since $p_g = q = 0$ for rational surface $X$, this comes
immediately from Theorem \ref{thm:observation}. For the last
statement, see Remark 1.4 in $\S 1$. \qed \\

\noindent {\bf Proof of Corollary \ref{cor:generalsubsystem}.} For
the complete embedding $X \hookrightarrow \P H^0 (X,L)$, assume
that $C$ is embedded in the linear subspace $\Lambda=\P H^0 (C,L
\otimes \mathcal{O}_C) \subset \P H^0 (X,L)$. If the center of the
linear projection $\P H^0 (X,L) \dashrightarrow \P (V)$ is $\Gamma
\cong \P ^{t-1} \subset \P H^0 (X,L)$, then the image of the map
$V \rightarrow H^0 (C,L \otimes \mathcal{O}_C)$ has codimension
$\leq m-2$ if and only if $\mbox{dim}(\Gamma \cap \Lambda) \leq
m-3$. For general $\Gamma \subset \P H^0 (X,L)$, this holds if and
only if
\begin{equation*}
\mbox{dim}(V) \geq \mbox{dim} (\Lambda) -m +3
\end{equation*}
or equivalently
\begin{equation*}
\mbox{codim}(V,H^0 (X,L)) \leq h^0 (X,L\otimes \mathcal{O}_X (-C))
+m-2
\end{equation*}
since $\mbox{dim}(\Lambda) = h^0 (C,L \otimes \mathcal{O}_C)-1 =
h^0(X,L) - h^0 (X,L\otimes \mathcal{O}_X (-C))$. This proves the
last statement of Theorem \ref{thm:observationrational}.\qed\\

\section{Anticanonical Rational Surfaces}
\noindent In this section we apply Theorem
\ref{thm:observationrational} to anticanonical rational surfaces.
Note that Theorem \ref{thm:observationrational} can be applied to
all very ample line bundles on anticanonical rational surfaces.

\subsection{Projective Plane} Let $Q \subset \P^2$
be a smooth quadric in $|\mathcal{O}_{\P^2} (2)|$. Let $d \geq 2$
be a natural number. By the Veronese embedding $\P^2
\hookrightarrow \P H^0 (\P^2,\mathcal{O}_{\P^2} (d))$, $Q$ is
contained in a linear subspace $\Lambda \subset \P H^0
(\P^2,\mathcal{O}_{\P^2} (d))$ of codimension ${{d} \choose {2}}$.

\begin{theorem}\label{thm:plane}
$(1)$ For a very ample subsystem $V \subset H^0
(\P^2,\mathcal{O}_{\P^2} (d))$, if the image of

the map $V \rightarrow H^0 (Q,\mathcal{O}_{\P^2} (d) \otimes
\mathcal{O}_Q)=H^0 (\P^1,\mathcal{O}_{\P^1} (2d))$ has codimension
$\leq 2$,

then $\P^2 \hookrightarrow \P(V)$ satisfies Property $N^S_1$.\\
$(2)$ General very ample subsystem $V \subset H^0
(\P^2,\mathcal{O}_{\P^2} (d))$ of dimension $\geq 2d$ satisfies

Property $N^S _1$.
\end{theorem}

\begin{proof}
The assertions follows immediately from Theorem
\ref{thm:observationrational} and Corollary \ref
{cor:generalsubsystem}.
\end{proof}

\noindent {\bf Remark 4.1} G. Ottaviani and R. Paoletti \cite{OP}
proved that for $d \geq 3$, $(\P^2,\mathcal{O}_{\P^2} (d))$
satisfies Property $N_p$ if and only if $d \leq 3d-3$. Thus if $V
\subset H^0 (\P^2,\mathcal{O}_{\P^2} (d))$ is a very ample
subsystem of codimension $\leq 3d-4$, then $\P^2 \hookrightarrow
\P (V)$ satisfies Property $N^S _1$ by Theorem \ref{thm:effect}.
\qed

\subsection{Hirzebruch Surfaces}
Throughout this subsection let $X = \F _e$, for any $e \geq 0$, be
the Hirzebruch surface defined by $\mathcal{O}_{\P^1} \oplus
\mathcal{O}_{\P^1} (-e)$. For the minimal section $C_0$ and a
fiber $f$, $C_0 + (e+1)f$ is a very ample line bundle and a
general section $C \in |C_0 + (e+1)f|$ is a smooth rational curve
with $C^2 = e+2$.

\begin{theorem}\label{thm:Hirzebruch}
Let $L=aC_0 +bf$ be a very ample line bundle. \\
$(1)$ For a very ample subsystem $V \subset H^0 (X,L)$, if the
image of the natural map

$V \rightarrow H^0 (C,L \otimes \mathcal{O}_C)$ has codimension
$\leq e$, then $X \hookrightarrow \P(V)$ satisfies Property $N^S_1$.\\
$(2)$ General very ample subsystem $V \subset H^0 (X,L)$ of
dimension $\geq a+b-e+1$ satisfies

Property $N^S _1$.
\end{theorem}

\begin{proof}
Note that for every very ample line bundle $L \in \mbox{Pic}X$,
$H^1 (X,L^j)=0$ for every $j \geq 1$ and thus the assertions
follow immediately from Theorem \ref{thm:observationrational},
Corollary \ref{cor:generalsubsystem} and Riemann-Roch.
\end{proof}

\noindent {\bf Remark 4.2.} F. J. Gallego and  B. P.
Purnaprajna\cite{GP} proved that $(\F_e,aC_0 +bf)$ with $e \geq 1$
and $a\geq2$ or $e=0$ and $a,b \geq 2$ satisfies Property $N_p$ if
and only if $2a+2b-ae \geq 3+p$. For these cases, if $V \subset
H^0 (\F_e,aC_0 +bf)$ is a very ample subsystem of codimension
$\leq 2a+2b-ae-4$, then the corresponding embedding satisfies
Property $N^S _1$.\qed

\subsection{Anticanonical Surfaces}
Let $X$ be an anticanonical rational surface with $K_X ^2 \leq 7$
and $L \in \mbox{Pic}X$ a very ample line bundle. Then $H^1
(X,L^j)=0$ for all $j \geq 1$. There is a birational morphism $X
\rightarrow \F_e$ for some $e \geq 0$. Let $C \subset X$ be the
preimage of a general section in $|C_0 + (e+1)f|$.

\begin{theorem}\label{thm:anticanonical}
$(1)$ Under the situation just stated, let $V \subset H^0 (X,L)$ a
very ample

subsystem. If the image of the natural map $V \rightarrow H^0 (C,L
\otimes \mathcal{O}_C)$ has codimension

$\leq e$, then $X \hookrightarrow \P(V)$ satisfies Property $N^S_1$.\\
$(2)$ Let $V \subset H^0 (X,L)$ be a very ample subsystem  such
that
\begin{equation*}
\mbox{dim} V\geq h^0 (C,L \otimes \mathcal{O}_C)-h^1 (X,L \otimes
\mathcal{O}_X (-C))-e
\end{equation*}

Then $X \hookrightarrow \P (V)$ satisfies Property $N^S _1$.
\end{theorem}

\begin{proof}
The assertions follows immediately from Theorem
\ref{thm:observationrational}, Corollary
\ref{cor:generalsubsystem} and Riemann-Roch.
\end{proof}

\noindent {\bf Remark 4.3.} F. J. Gallego and  B. P.
Purnaprajna\cite{GP} proved that $(X,L)$  satisfies Property $N_p$
if and only if $-L \cdot K_X \geq 3+p$. For these cases, if $V
\subset H^0 (X,L)$ is a very ample subsystem of codimension $\leq
p-1$, then the corresponding embedding satisfies Property $N^S
_1$.\qed

\end{document}